\documentclass[12pt,a4paper]{amsart}
\usepackage{amssymb}
\usepackage{amscd}
\begin{document}

\setcounter{MaxMatrixCols}{30}
\newtheorem{lemma}{Lemma}[section]
\newtheorem{prop}[lemma]{Proposition}
\newtheorem{cor}[lemma]{Corollary}
\newtheorem{thm}[lemma]{Theorem}
\newtheorem{con}[lemma]{Conjecture}
\theoremstyle{definition}
\newtheorem{rem}[lemma]{Remark}
\newtheorem{rems}[lemma]{Remarks}
\newtheorem{defi}[lemma]{Definition}
\newtheorem{ex}[lemma]{Example}
\title[Bass conjecture for amenable groups]{From acyclic groups\\ to\\ the Bass conjecture for amenable groups}
\author{A. J. Berrick, I. Chatterji and G. Mislin}
\address{\textsc{Department of Mathematics, National University of Singapore,
Singapore\qquad\qquad\qquad\qquad\qquad\qquad\qquad\qquad\qquad\qquad
\qquad\qquad\hfill}\qquad\qquad\qquad\\
\textsc{Department of Mathematics, Cornell University, White Hall,
Ithaca NY 14853-7901,
USA\qquad\qquad\qquad\qquad\qquad\qquad\qquad\qquad\qquad\qquad
\qquad\qquad\hfill}\qquad\qquad\qquad\\
\textsc{Department of Mathematics, ETH Z\"{u}rich, 8092
Z\"{u}rich,
Switzerland}\\
\texttt{berrick@math.nus.edu.sg\quad indira@math.cornell.edu\quad
mislin@math.ethz.ch}}
\date{9 January, 2004}

\begin{abstract}
We prove that the Bost Conjecture on the $\ell^{1}$-assembly map for countable
discrete groups implies the Bass Conjecture. It follows that all amenable
groups satisfy the Bass Conjecture.

\end{abstract}
\maketitle

\section{Introduction}

Throughout this paper, let $G$ be a discrete group. For each finitely
generated projective (left) $\mathbb{Z}G$-module $P$, there exists an
idempotent matrix $(m_{ij})=M\in M_{n}(\mathbb{Z}G)$ such that $P$ is
isomorphic to the image under right multiplication $\mathbb{Z}G^{n}%
\rightarrow\mathbb{Z}G^{n}$ by $M$. Writing $[\mathbb{Z}G,\mathbb{Z}G]$ for
the additive subgroup of $\mathbb{Z}G$ generated by the elements $gh-hg$
($g,h\in G$), we identify $\mathbb{Z}G/[\mathbb{Z}G,\mathbb{Z}G]$ with
$\bigoplus_{[s]\in\lbrack G]}\mathbb{Z}\cdot\lbrack s]$, where $[G]$ is the
set of conjugacy classes of elements of $G$. The \emph{Hattori-Stallings rank}
$r_{P}$ is then defined by
\begin{align*}
r_{P} =\sum_{i=1}^{n}m_{ii}+[\mathbb{Z}G,\mathbb{Z}G]
=\sum_{\lbrack s]\in\lbrack G]}r_{P}(s)[s]\in\bigoplus_{\lbrack
s]\in\lbrack G]}\mathbb{Z}\cdot\lbrack s].
\end{align*}
In his seminal paper \cite{bass}, H. Bass made the following conjecture.

\begin{con}[Classical Bass Conjecture] For any finitely generated
projective $\mathbb{Z}G$-module $P$, the values $r_{P}(s)\in\mathbb{Z}$ of the
Hattori-Stallings rank $r_{P}$ are zero for $s\in G\setminus\{1\}$.
\end{con}

One of the striking results of \cite{bass} is a proof of this conjecture for
the case of torsion-free linear groups (for the case of arbitrary linear
groups see \cite{Linnell}). Using methods of cyclic homology, Eckmann in
\cite{Eckmann1} and Emmanouil in \cite{Emmanouil} proved the conjecture for
many more groups (for a recent survey see Eckmann \cite{Eckmann}). The case of
solvable groups was settled only very recently by Farrell and Linnell
\cite{fl}, where they prove the classical Bass conjecture for all elementary
amenable groups.

We show below that groups that satisfy the \emph{Bost conjecture} satisfy the
classical Bass conjecture too, and indeed more general versions thereof that
we call the $\ell^{1}$ Bass conjecture~(\ref{L1Bass}) and $\mathbb{C}G$ Bass
conjecture. Combining with known information concerning the Bost conjecture
gives the result of the title.

\begin{thm}
\label{pain}Amenable groups satisfy the classical Bass conjecture.
\end{thm}

The class of amenable groups includes the class of elementary
amen\-able groups, which is the class obtained from finite and
abelian groups by means of subgroups, quotients, extensions and
increasing unions. This inclusion is strict \cite{Grigorchuk}.

The proof is obtained via the following chain of deductions, describing a tour
from geometric functional analysis, through operator algebra $K$-theory,
algebraic topology and combinatorial group theory, and ultimately to general
linear algebra.

\begin{thm}[Lafforgue \cite{Lafforgue}]
\label{Lafforgue: Bost} For any countable discrete group with
the Haagerup property (for example, any countable amenable group
\cite{BeCheVal}), the Bost assembly map
\[
\beta_{\ast}^{G}:K_{\ast}^{G}(\underline{E}G)\rightarrow K_{\ast
}^{\mathrm{top}}(\ell^{1}(G))
\]
is an isomorphism.
\end{thm}

We refer to Lafforgue's work in \cite{Lafforgue} for the definition of the
Bost assembly map. The Bost conjecture (see \cite{Skandalis}) is that the
isomorphism holds for all countable discrete groups.

The brunt of our paper consists in proving the following key link between
these conjectures.

\begin{thm}
\label{Bost=>L1Bass copy(1)} Let $G$ be a countable discrete group for which
the Bost assembly map is rationally an epimorphism in degree $0$. Then the
$\ell^{1}$ Bass conjecture \ref{L1Bass} holds for $G$.
\end{thm}

The proof of this theorem involves a natural embedding of $G$ in an acyclic
group $A(G)$ that is injective on conjugacy classes, with the centralizer of
any finitely generated abelian subgroup of $A(G)$ acyclic as well (cf.
Theorem~\ref{Berrick} below). This allows us to control the image of the
universal trace
\[
T^{1}:K_{0}(\ell^{1}(G))\rightarrow HH_{0}(\ell^{1}(G)).
\]

The proof of the classical Bass conjecture is then clinched by an easy lemma.

\begin{lemma}
\label{BassCjs}\emph{(a) }\label{countable copy(1)}The $\ell^{1}$ Bass
conjecture \ref{L1Bass} holds for a group $G$ if it holds for all its
countable subgroups.

\begin{enumerate}
\item[\textrm{(b)}] \label{L1Bass => CG Bass}The $\ell^{1}$ Bass conjecture
\ref{L1Bass} implies the $\mathbb{C}G$ Bass conjecture \ref{Bass}.

\item[\textrm{(c)}] \label{CG Bass => classical Bass}The $\mathbb{C}G$ Bass
conjecture \ref{Bass} implies the classical Bass conjecture.
\end{enumerate}
\end{lemma}

Theorem \ref{Bost=>L1Bass copy(1)} also provides information concerning the
non-existence of idempotents in $\mathbb{C}G$ (Kaplansky conjecture), and
indeed an analogous result for the case of $\ell^{1}(G)$. Of course, here one
needs to assume the group in question to be torsion-free.

\begin{cor}
\label{B} Let $G$ be a torsion-free countable group and assume that the Bost
assembly map $\beta_{0}^{G}$ is rationally surjective. Then $\ell^{1}(G)$
contains no idempotent other than $0$ and $1$.
\end{cor}

For $G$ a torsion-free hyperbolic group, Ji proved already in
\cite{Ji} that $\ell^{1}(G)$ does not have any idempotent other
than $0$ and $1$. The analogous result for $C_{r}^{\ast}(G)$ was
recently proved by Puschnigg \cite{Pu} as well as Mineyev and Yu
\cite{MiYu}. As pointed out to us by A. Valette, this corollary
can easily be deduced using the Atiyah $L^{2}$-index theorem, but
our proof avoids it, and as a matter of fact our techniques allow
a new proof for this theorem, see \cite{ell2}. As an application
of our methods, we obtain a new proof of a recent result of
L\"{u}ck \cite{Lueck}, which is a bound for the image of the
composite of the Kaplansky trace following the Baum-Connes
assembly map (see Section \ref{ImageTrace}).

For a discussion of other groups to which our argument applies,
see Section \ref{ClasseC'} below. As well as groups that have the
Haagerup property, they include discrete subgroups of virtually
connected semisimple linear Lie groups, hyperbolic groups and
cocompact CAT(0)-groups, to name a few.


\section{Traces on $\mathbb{C}G$ and $\ell^{1}(G)$}

Let $A$ be an algebra over $\mathbb{C}$ and $M$ be a $\mathbb{C}$-module. A
\emph{trace map} over $A$ is a $\mathbb{C}$-linear map $t:A\rightarrow M$
satisfying $t(ab)=t(ba)$ for any $a,b\in A$. There is a \emph{universal
trace}
\[
T:A\rightarrow HH_{0}(A)=A/[A,A],\ a\mapsto a+[A,A]
\]
where $HH_{0}(A)$ stands for the $0$-th Hochschild homology group, and $[A,A]$
denotes the $\mathbb{C}$-submodule of $A$, generated by the elements of the
form $ab-ba$ (for $a,b\in A$). In fact, every trace factors uniquely through
$T$. From consideration of idempotent matrix representatives as in the
Introduction, any trace gives rise to a well-defined map $K_{0}(A)\rightarrow
M$, where $K_{0}(A)$ is the Grothendieck group of finitely generated
projective $A$-modules. In particular, the universal trace induces
\[
T:K_{0}(A)\rightarrow HH_{0}(A).
\]
We focus on the cases where $A$ is the group algebra $\mathbb{C}G$ or the
Banach algebra $\ell^{1}(G)$ of summable series $a=\sum_{g\in G}a_{g}g$, with
$\ell^{1}$-norm $\Vert a\Vert_{1}=\sum_{g\in G}|a_{g}|<\infty$, for $G$ a
discrete group. We sometimes refer to elements $a=\sum_{g\in G}a_{g}g\in
\ell^{1}(G)$ as functions $G\rightarrow\mathbb{C}$, $g\mapsto a_{g}$. In this
way $\mathbb{C}G\subset\ell^{1}(G)$ corresponds to the functions with finite
support. For $A=\mathbb{C}G$ or $\ell^{1}(G)$ we now define the
\emph{Kaplansky trace}
\begin{align*}
\kappa:A  &  \rightarrow\mathbb{C}\\
a  &  \mapsto a_{1}%
\end{align*}
if $a=\sum_{g\in G}a_{g}g$ and $1$ denotes the neutral element in $G$; note
that indeed $\kappa(ab)=\kappa(ba)$. We use the same notation for the induced
map $\kappa:K_{0}(A)\rightarrow\mathbb{C}$. Another example of trace in these
cases is given by the \emph{augmentation trace}
\begin{align*}
\epsilon:A  &  \rightarrow\mathbb{C}\\
a  &  \mapsto\sum_{g\in G}a_{g}.
\end{align*}
The special feature of that trace is that it is an algebra homomorphism, and
thus, $K_{0}(\ )$ and $HH_{0}(\ )$ being covariant functors from the category
of algebras to the category of abelian groups, we get the following
commutative diagram:
\[
\begin{CD}
K_0(A)               @>{T}>>    HH_0(A)\\
@V{K_0(\epsilon)}VV             @VV{HH_0(\epsilon)}V\\
K_0(\mathbb{C})              @>{T}>>    HH_0(\mathbb{C})=\mathbb{C}.
\end{CD}\label{T from A to C}%
\]
Since $K_{0}(\mathbb{C})=\mathbb{Z}$, we deduce that the induced map
$\epsilon:K_{0}(A)\rightarrow\mathbb{C}$ has image $\mathbb{Z}\subset
\mathbb{C}$.

We consider a third trace, which we discuss first for $A=\mathbb{C}G$. We call
as usual \emph{Hattori-Stallings trace} the map
\begin{align*}
HS:\mathbb{C}G  &  \rightarrow\bigoplus_{[G]}\mathbb{C}\\
a  &  \mapsto\sum_{[x]\in\lbrack G]}\epsilon_{\lbrack x]}(a)[x]
\end{align*}
where $[x]$ stands for the conjugacy class of $x\in G$, $[G]$ the set of
conjugacy classes of $G$ and $\epsilon_{\lbrack x]}(a)=\sum_{g\in\lbrack
x]}a_{g}$. This is in fact nothing else but the universal trace for the case
$A=\mathbb{C}G$. The $[x]$-component $HS_{[x]}$ of the associated trace
\[
HS:K_{0}(\mathbb{C}G)\rightarrow\bigoplus_{\lbrack G]}\mathbb{C}%
\]
satisfies $HS_{[x]}([\mathbb{C}G\otimes_{G}P])=r_{P}(x)$, where $P$ denotes a
finitely generated projective $\mathbb{Z}G$-module. We recall that, by a
result of Linnell (see \cite{Linnell}, Lemma 4.1), one has $r_{P}(x)=0$ for
$x\in G\setminus\{1\}$ of finite order. Thus we consider the following version
of the Bass conjecture.

\begin{con}
[$\mathbb{C}G$ Bass Conjecture]\label{Bass} The
Hattori-Stallings trace
\[
HS:K_{0}(\mathbb{C}G)\rightarrow\bigoplus_{\lbrack G]}\mathbb{C}%
\]
takes its values in the $\mathbb{C}$-vector space spanned by the conjugacy
classes of elements of finite order.
\end{con}

\noindent\textit{Proof of Lemma \ref{CG Bass => classical Bass}(c). }To deduce
from the above conjecture the classical Bass conjecture concerning
$\mathbb{Z}G$, one considers the natural map
\[
i_{\ast}:K_{0}(\mathbb{Z}G)\rightarrow K_{0}(\mathbb{C}G)
\]
and it suffices (because of Linnell's result) to show that, for $[P]\in
K_{0}(\mathbb{Z}G)$ and for $x\in G$ of infinite order, $HS_{[x]}(i_{\ast
}[P])=0$. Therefore Conjecture \ref{Bass} implies the classical Bass
conjecture concerning $\mathbb{Z}G$.\hfill$\square$\medskip

In the case where $A=\ell^{1}(G)$, per definition, the universal trace $T$ on
$K_{0}(\ell^{1}(G))$ takes its values in $HH_{0}(\ell^{1}(G))$; the topology
on $\ell^{1}(G)$ does not enter here. Now consider the Banach space $\ell
^{1}([G])$, the completion of the vector space $\bigoplus_{[G]}\mathbb{C}$
with respect to the $\ell^{1}$-norm. Again, we think of elements $\sum
_{[x]\in\lbrack G]}a_{[x]}[x]$ in $\ell^{1}([G])$ as functions $[G]\rightarrow
\mathbb{C}$, $[x]\mapsto a_{[x]}$. If we write $\mathrm{FC}(G)\subset\lbrack
G]$ for the subset of conjugacy classes $[g]\in\lbrack G]$ with $g$ of finite
order, we have
\[
\bigoplus_{\mathrm{FC}(G)}\mathbb{C}\subset\bigoplus_{\lbrack G]}%
\mathbb{C}\subset\ell^{1}([G]).
\]
For $a=\sum_{g\in G}a_{g}g\in\ell^{1}(G)$, we define $\epsilon_{\lbrack
x]}(a)\in\mathbb{C}$ by $\epsilon_{\lbrack x]}(a)=\sum_{g\in\lbrack x]}a_{g}$,
which is a convergent series because $\sum_{g\in G}|a_{g}|<\infty$. The map
\begin{align*}
p:\ell^{1}(G)  &  \rightarrow{\ell^{1}([G])}\\
a  &  \mapsto\sum_{\lbrack x]\in\lbrack G]}\epsilon_{\lbrack x]}(a)[x]
\end{align*}
is a trace and therefore induces a well-defined map
\[
\overline{p}:HH_{0}(\ell^{1}(G))\rightarrow{\ell^{1}([G])}.
\]
Indeed, any $z\in\lbrack\ell^{1}(G),\ell^{1}(G)]$ has the form%
\[
z=\sum_{i=1}^{n}(a^{i}b^{i}-b^{i}a^{i})=\sum_{i=1}^{n}\sum_{g,h\in G}a_{g}%
^{i}b_{h}^{i}(gh-hg)
\]
for $a^{i}=\sum_{g\in G}a_{g}^{i}g$ and $b^{i}=\sum_{h\in G}b_{h}^{i}h$
elements in $\ell^{1}(G)$; thus, since $[gh]=[hghh^{-1}]=[hg],$ for all $x\in
G$ the $[x]$-component of $p(z)$ is given by the absolutely convergent series
$\sum_{i=1}^{n}\sum_{gh\in\lbrack x]}a_{g}^{i}b_{h}^{i}-a_{g}^{i}b_{h}^{i}=0$.
Hence the trace $p$ yields an $\ell^{1}$-version of the Hattori-Stallings
trace
\[
HS^{1}:=\overline{p}\circ T:K_{0}(\ell^{1}(G))\rightarrow{\ell^{1}([G])}%
\]
and leads us to the following $\ell^{1}$-version of the Bass conjecture.

\begin{con}[$\ell^{1}$ Bass Conjecture]
\label{L1Bass}The $\ell^{1}$
Hattori-Stallings trace
\[
HS^{1}:K_{0}(\ell^{1}(G))\rightarrow{\ell^{1}([G])}%
\]
takes its values in the subspace $\bigoplus_{\mathrm{FC}(G)}\mathbb{C}$ of
functions finitely supported by the conjugacy classes of elements of finite order.
\end{con}

\noindent\textit{Proof of Lemma \ref{L1Bass => CG Bass}(b). }With
$j:\mathbb{C}G\hookrightarrow\ell^{1}(G)$, the commutative diagram
\[
\begin{CD}
K_0(\mathbb{C}G)  @>HS>> {\bigoplus\limits_{[G]}\mathbb{C}}\\
@V{j}_*VV      @VVV\\
K_0(\ell^1(G))@>HS^1>>{\ell^1([G])}
\end{CD}\label{HS from C to l1}%
\]
shows that if $HS^{1}$ maps into $\bigoplus_{\mathrm{FC}(G)}\mathbb{C}$, then
so does $HS$.\hfill$\square$\medskip

Certain arguments in the sequel require $G$ to be a countable group. For our
applications to the Bass conjecture, this does not create any problem, in view
of the following.\medskip

\noindent\textit{Proof of Lemma \ref{countable copy(1)}(a). }Notice that an
idempotent matrix in $M_{n}(\ell^{1}(G))$ representing a finitely generated
projective $\ell^{1}(G)$-module $P$ involves only countably many elements from
$G$. So, for some countable subgroup $G_{\alpha}$ of $G$ and finitely
generated projective ${\ell^{1}(G_{\alpha})}$-module $Q$, $P$ is of the form
$\ell^{1}(G)\otimes_{\ell^{1}(G_{\alpha})}Q$. Thus the inclusion maps
$G_{\alpha}\hookrightarrow G$ of all countable subgroups $G_{\alpha}$ of $G$
induce an epimorphism
\[
\bigoplus\limits_{\alpha}K_{0}(\ell^{1}(G_{\alpha}))\rightarrow K_{0}(\ell
^{1}(G)),
\]
and naturality of $HS^{1}$ finishes the argument.\hfill$\square$\medskip

Of course, although not needed here, a similar argument also holds for
Conjecture \ref{Bass}.


\section{$K^{G}_{0}\Lambda$-discrete spaces}

We recall some more terminology (discussed further in, for example,
\cite{Mislin}).

The \emph{reduced C*-algebra} $C^{*}_{r}(G)$ is the completion of $\ell
^{1}(G)$ with respect to the operator norm, where $\ell^{1}(G)$ acts on the
Hilbert space $\ell^{2}(G)$ of square summable functions via its regular
representation. Recall that $C^{*}_{r}(\ )$ is not a functor on the category
of groups, but for $H<G$ a subgroup, $C^{*}_{r}(H)$ is in a natural way a
subalgebra of $C^{*}_{r}(G)$, so that for an injective group homomorphism one
does obtain an induced morphism of reduced $C^{*}$-algebras.

\medskip

Let $KK_{\ast}^{G}(A,B)$ denote the equivariant Kasparov $K$-groups of the
pair of separable $G$-$C^{\ast}$-algebras $A,B$ (see \cite{Kasparov}). Recall
that a $G$-CW-complex $X$ is said to be \emph{proper} if the stabilizer of
each vertex is finite. For $X$ a proper $G$-CW-complex with $G$ a
\emph{countable} group, the equivariant $K$-homology groups $RK_{\ast}^{G}(X)$
are then defined by
\[
RK_{\ast}^{G}(X)=\text{\textrm{colim}}_{\{Y\subset X\mid Y\,{\text{cocompact}%
}\}}\,KK_{\ast}^{G}(C_{0}Y,\,\mathbb{C}),
\]
where $Y$ runs over the cocompact $G$-subcomplexes of $X$, $\mathbb{C}$ is
considered as a $C^{\ast}$-algebra with trivial $G$-action, and $C_{0}Y$
denotes the $C^{\ast}$-algebra of continuous functions $Y\rightarrow
\mathbb{C}$ on the locally compact CW-complex $Y$ that vanish at infinity.

\smallskip

These homology groups $RK_{*}^{G}$ turn out to be representable in
the following sense. We write $\mathfrak{O}(G)$ for the orbit
category of $G$ (the objects are the cosets $G/H$ and morphisms
are $G$-maps). There exists an $\mathfrak{O}(G)$-spectrum
representing a homology theory $K_{*}^{G}$ on the category of
\emph{all} $G$-CW-complexes such that for $H<G$ one has 
\[
K_{*}^{G}(G/H)=K_{*}^{\mathrm{top}}(C^{*}_{r}H),
\]
where $K_{*}^{\mathrm{top}}(C^{*}_{r}H)$ is the (topological) algebraic
$K$-theory of the Banach algebra $C^{*}_{r}H$ (note that $K_{0}^{\mathrm{top}%
}(C^{*}_{r}G)=K_{0}(C^{*}_{r}G)$, the projective class group of the ring
$C^{*}_{r}G$ -- the topology does not enter in this case, see \cite{Karoubi}).
Moreover, the homology theory $K_{*}^{G}$ is such that for all \emph{proper}
$G$-CW-complexes $X$ and countable groups $G$ one has a natural isomorphism
\[
K_{*}^{G}(X)\cong RK_{*}^{G}(X).
\]
We emphasize that the right-hand side is defined only in case $G$ is
countable, whereas the left-hand side is defined for any discrete group $G$.
Note that, because $K^{G}_{*}$ is defined by a spectrum, it is fully
additive:
\[
K^{G}_{*}(\coprod X_{\alpha})\simeq\bigoplus K^{G}_{*}(X_{\alpha}).
\]
For details concerning the $\mathfrak{O}(G)$-spectrum representing equivariant
$K$-homology, the reader is referred to Davis and L\"uck \cite{davislueck}.

\begin{defi}
Let $\Lambda$ be a (unital) subring of $\mathbb{C}$. A proper $G$-CW-complex
$X$ is called \emph{$K^{G}_{0}\Lambda$-discrete} if the natural map induced by
the inclusion $\iota:X^{0}\hookrightarrow X$ of the $0$-skeleton
\[
\iota_{*}:K^{G}_{0}(X^{0})\to K^{G}_{0}(X)
\]
is an epimorphism after tensoring with $\Lambda$.
\end{defi}

Clearly, being $K^{G}_{0}\Lambda$-discrete depends only on the $G$-homotopy
type of the $G$-CW-complex $X$.

\begin{rem}
\label{0-sk}If $X$ is a proper $G$-CW-complex, then $X^{0}=\coprod
G/F_{\alpha}$ for some set of finite subgroups $F_{\alpha}\subset G$.
Therefore
\[
K_{0}^{G}(X^{0})=K_{0}^{G}({\small {\coprod}}G/F_{\alpha})=\bigoplus K_{0}%
^{G}(G/F_{\alpha})
\]
and, since $G/F_{\alpha}$ is the $G$-space induced from the proper $F_{\alpha
}$-space $\{\mathrm{pt}\}$,
\[
K_{0}^{G}(G/F_{\alpha})\cong K_{0}^{F_{\alpha}}(\{\mathrm{pt}\})\cong
R_{\mathbb{C}}(F_{\alpha}),
\]
where $R_{\mathbb{C}}(F_{\alpha})$ is the complex representation ring of the
finite group $F_{\alpha}$. Thus one has a natural isomorphism of abelian
groups
\[
\beta:K_{0}^{G}(X^{0})\rightarrow\bigoplus R_{\mathbb{C}}(F_{\alpha}).
\]

\end{rem}

Recall that a group $G$ is called \emph{acyclic} if the classifying space
$K(G,1)$ $=$ $BG$ satisfies $H_{*}(BG;\mathbb{Z})=0$ for $*>0$. Equivalently,
$G$ is acyclic if the suspension $\Sigma BG$ of $BG$ is contractible. As
usual, we write $EG$ for the universal cover of $BG$. It is a free
$G$-CW-complex, and so proper.

\begin{lemma}
If $G$ is acyclic, then $EG$ is $K^{G}_{0}\mathbb{Z}$-discrete.
\end{lemma}

\begin{proof}
Using a suitable model for $EG$ we may assume that $EG^{0}=G$ as a discrete
$G$-space, and therefore
\[
K^{G}_{0}(EG^{0})=K^{G}_{0}(G)=K_{0}(\{\mathrm{pt}\})\cong\mathbb{Z}.
\]
On the other hand, $K^{G}_{0}(EG)\cong K_{0}(BG)$ and, as here the suspension
$\Sigma BG$ is contractible, the inclusion $\{\mathrm{pt}\}\hookrightarrow BG$
induces an isomorphism
\[
K^{G}_{0}(EG^{0})\cong K_{0}(\{\mathrm{pt}\})\cong K_{0}(BG)\cong K^{G}%
_{0}(EG),
\]
showing that $EG$ is $K^{G}_{0}\mathbb{Z}$-discrete.
\end{proof}

The universal proper $G$-CW-complex $\underline{E}G$ is characterized up to
$G$-homotopy by the property that
\[
(\underline{E}G)^{H}\simeq\left\{
\begin{array}
[c]{cc}%
\{\mathrm{pt}\} & \text{if }|H|<\infty\\
\varnothing & \text{otherwise.}%
\end{array}
\right.
\]
This implies that for any finite subgroup $H<G$, with centralizer $C_{G}%
(H)<G$, the $C_{G}(H)$-CW-complex $(\underline{E}G)^{H}$ is a model for
$\underline{E}C_{G}(H)$. For a discussion of $\underline{E}G$ and its
properties, see for instance \cite{PeterGuido}. If $G$ is torsion-free, then
$EG$ is a model for $\underline{E}G$, so that $\underline{E}G$ is $K_{0}%
^{G}\mathbb{Z}$-discrete when $G$ is also acyclic. In the context
of the Bost conjecture or Baum-Connes conjecture (\ref{BC} below),
we are interested in $\underline{E}G$, which differs from $EG$ as
soon as $G$ is not torsion-free. So, to deal with groups with
torsion, we need a stronger version of acyclicity, which in
particular takes into account the centralizers of finite order
elements in the group.

\begin{defi}
\label{pervasively} A group $G$ is called \emph{pervasively acyclic}, if for
every finitely generated abelian subgroup $A<G$ the centralizer $C_{G}(A)$ is
an acyclic group.
\end{defi}

We also need a way of keeping track of the torsion in $G$.

\begin{defi}
For any group $G$, $\Lambda_{G}$ denotes the subring of $\mathbb{Q}$ generated
by the elements $1/|H|$, where $H$ runs over the finite subgroups of $G$.
\end{defi}

The following lemma is useful later.

\begin{lemma}
\label{BrownSpectral}For a group $G$, the $G$-map $EG\to\underline{E}G$
induces an isomorphism
\[
H_{*}(EG/G;\Lambda)\to H_{*}(\underline{E}G/G;\Lambda)
\]
for any abelian group $\Lambda$ such that $\Lambda_{G}\subset\Lambda
\subset\mathbb{Q}$.
\end{lemma}

\begin{proof}
The Brown spectral sequence (see \cite{Brown}, VII, 7.10) for the
$G$-CW-complex $\underline{E}G$ takes the form
\[
E^{1}_{p,q}=\bigoplus_{\sigma\in\Sigma_{p}}H_{q}(BG_{\sigma};\Lambda
)\Longrightarrow H_{p+q}(BG;\Lambda)
\]
with $G_{\sigma}$ the (finite) stabilizer of the $p$-cell $\sigma
\subset\underline{E}G$ and $\Sigma_{p}$ a set of representatives of orbits of
$p$-cells. Our assumption on $\Lambda$ implies that $E^{1}_{p,q}=0$ for $q>0$,
so that $E^{1}_{p,*}=E^{1}_{p,0}\cong C_{p}(\underline{E}G/G;\Lambda)$, the
cellular chain complex of $\underline{E}G/G$ with coefficients in $\Lambda$.
Thus
\[
E^{2}_{p,*}=E^{2}_{p,0}\cong H_{p}(\underline{E}G/G;\Lambda)\cong E^{\infty
}_{p,0}\cong H_{p}(BG;\Lambda),
\]
with the isomorphism being induced by the edge-homomorphism.
\end{proof}

Our aim now is to show that, for any pervasively acyclic group $G$, its
classifying space for proper actions $\underline{E}G$ is $K^{G}_{0}\Lambda
_{G}$-discrete.

We first recall the use of Bredon homology in the context of equivariant
$K$-homology (we refer to \cite{Mislin} for a more detailed exposition of
these techniques). Let $G$ be any group and $\mathfrak{F}$ the set of finite
subgroups of $G$. We write $\mathfrak{O}(G,\mathfrak{F})$ for the full
subcategory of the orbit category $\mathfrak{O}(G)$ with objects $G/H$,
$H\in\mathfrak{F}$. If $X$ is a $G$-CW-complex, its cellular chain complex
$C_{\ast}X$ gives rise to a contravariant functor into the category
$\mathrm{Ab}_{\ast}$ of chain complexes of abelian groups,%
\[
C_{\ast}^{\mathfrak{F}}X:\mathfrak{O}(G,\mathfrak{F})\longrightarrow
\mathrm{Ab}_{\ast},\qquad G/H\longmapsto C_{\ast}X^{H}.
\]
If $M:\mathfrak{O}(G,\mathfrak{F})\longrightarrow\mathrm{Ab}$ is a (covariant)
functor with target $\mathrm{Ab}$ (the category of abelian groups), the Bredon
homology $H_{\ast}^{\mathfrak{F}}(X;M)$ is defined to be the homology of the
chain complex of abelian groups $C_{\ast}^{\mathfrak{F}}X\otimes
_{\mathfrak{F}}M.$ The latter is defined as
\[
\sum\limits_{H\in\mathfrak{F}}\left.  C_{\ast}X^{H}\otimes M(G/H)\right/
\sim\,,
\]
where $\sim$ is the equivalence relation induced by $f^{\ast}x\otimes y\sim
x\otimes f_{\ast}y$, for $f:G/H\rightarrow G/K$ running over all morphisms of
$\mathfrak{O}(G,\mathfrak{F}).$

\begin{rem}
In the case of equivariant $K$-homology, the functor $M=K_{j}^{G}(?)$ is of
particular interest. For $j$ even, $K_{j}^{G}(G/H)\cong R_{\mathbb{C}}(H)$,
the complex representation ring of the finite group $H$, and for $j$ odd
$K_{j}^{G}(G/H)=0$. A morphism $\varphi:G/H\rightarrow G/K$ gives rise to
\[
\varphi_{\ast}:K_{0}^{G}(G/H)\rightarrow K_{0}^{G}(G/K),
\]
which corresponds to the map $R_{\mathbb{C}}(H)\rightarrow R_{\mathbb{C}}(K)$
induced by $h\mapsto x^{-1}hx$, if $\varphi(H)=xK$. In this way the
automorphism group
\[
\operatorname{map}_{G}(G/H,G/H)\cong N_{G}(H)/H
\]
acts on $K_{0}^{G}(G/H)=R_{\mathbb{C}}(H)$; in particular $C_{G}(H)$ acts
trivially on $R_{\mathbb{C}}(H)$ via $C_{G}(H)\rightarrow N_{G}(H)/H$. If $X$
is a proper $G$-CW-complex, there is an Atiyah-Hirzebruch spectral sequence
(cf. L\"{u}ck's Remark 3.9 in \cite{lueck2})
\[
E_{i,j}^{2}=H_{i}^{\mathfrak{F}}(X;\,K_{j}^{G}(?))\Longrightarrow K_{i+j}%
^{G}(X).
\]

\end{rem}

\begin{thm}
Let $X$ be a proper $G$-CW-complex. If for every finite cyclic subgroup $C<G$
one has $H_{2i}(X^{C}/C_{G}(C);\Lambda_{G})=0$ for all $i>0$, then $X$ is
$K_{0}^{G}\Lambda_{G}$-discrete.
\end{thm}

\begin{proof}
As in L\"{u}ck \cite{lueck2} we write ${\mathfrak{Sub}}(G,{\mathfrak{F}})$ for
the quotient category of $\mathfrak{O}(G,\mathfrak{F})$ with the same objects
as $\mathfrak{O}(G,\mathfrak{F})$, but with morphisms from $G/H$ to $G/K$
given by $\left.  \operatorname{map}_{G}(G/H,G/K)\right/  C_{G}(H)$, where the
centralizer $C_{G}(H)$ acts via
\[
C_{G}(H)\rightarrow N_{G}(H)/H\cong\operatorname{map}_{G}(G/H,G/H)
\]
on $G$-maps $G/H\rightarrow G/K$. The cellular chain complex of $X$ gives rise
to a contravariant functor to the category of chain complexes of abelian
groups:
\[
C_{\ast}^{\mathfrak{Sub}}X:{\mathfrak{Sub}}(G,{\mathfrak{F}})\rightarrow
\operatorname{Ab}_{\ast},\quad G/H\mapsto C_{\ast}(X^{H}/C_{G}(H)).
\]
As before, if $M:{\mathfrak{Sub}}(G,{\mathfrak{F}})\rightarrow
\operatorname{Ab}$ is any (covariant) functor, Bredon-type homology groups
$H_{\ast}^{\mathfrak{Sub}}(X;M)$ are defined. Note that the projection
$\pi:\mathfrak{O}(G,\mathfrak{F})\rightarrow{\mathfrak{Sub}}(G,{\mathfrak{F}%
})$ induces a functor $\pi^{\ast}M:=M\circ\pi$. There is a natural isomorphism
(cf. \cite{lueck2} (3.6))
\[
H_{\ast}^{\mathfrak{F}}(X;\,\pi^{\ast}M)\cong H_{\ast}^{\mathfrak{Sub}}(X;M).
\]
For the representation ring functor $R_{\mathbb{C}}:\mathfrak{O}%
(G,\mathfrak{F})\rightarrow\operatorname{Ab}$, $G/H\mapsto R_{\mathbb{C}%
}(H)=K_{0}^{G}(G/H)$, the centralizer $C_{G}(H)$ acts trivially on
$R_{\mathbb{C}}(H)$. Thus $R_{\mathbb{C}}$ factors through ${\mathfrak{Sub}%
}(G,{\mathfrak{F}})$, to yield a functor still denoted by $R_{\mathbb{C}}$,
and hence we can replace the left-hand side of the Atiyah-Hirzebruch spectral
sequence by $H_{i}^{\mathfrak{Sub}}(X;K_{j}^{G}(?))$. Denote by $\Lambda
_{G}{\mathfrak{S}}$ the category of functors ${\mathfrak{Sub}}(G,{\mathfrak{F}%
})\rightarrow\Lambda_{G}$-$\operatorname{Mod}$; an object of $\Lambda
_{G}{\mathfrak{S}}$ is called a $\Lambda_{G}{\mathfrak{S}}$-\emph{module}.
According to L\"{u}ck \cite{Lueck}, Theorem 3.5 (b), $\Lambda_{G}\otimes
R_{\mathbb{C}}$ is a projective $\Lambda_{G}{\mathfrak{S}}$-module. This
implies that for any proper $G$-CW-complex $X$
\[
H_{\ast}^{\mathfrak{Sub}}(X;\Lambda_{G}\otimes R_{\mathbb{C}})\cong H_{\ast
}^{\Lambda_{G}{\mathfrak{S}}}(X)\otimes_{\mathfrak{Sub}}(\Lambda_{G}\otimes
R_{\mathbb{C}}),
\]
where $H_{\ast}^{{\Lambda_{G}{\mathfrak{S}}}}(X)$ denotes the $\Lambda
_{G}{\mathfrak{S}}$-module $G/H\mapsto H_{\ast}(X^{H}/C_{G}(H);\Lambda_{G})$
(the `tensor product'\ $-\otimes_{\mathfrak{Sub}}-$ is, as before, defined by
taking the sum over all objects
\[
\sum_{H\in{\mathfrak{F}}}H_{\ast}(X^{H}/C_{G}(H);\Lambda_{G})\otimes\left(
\Lambda_{G}\otimes R_{\mathbb{C}}(H)\right)
\]
and dividing out by the equivalence relation generated by $f^{\ast}x\otimes
y\sim x\otimes f_{\ast}y$, with $f$ running over the morphisms of
${\mathfrak{Sub}}(G,{\mathfrak{F}})$).

According to Artin's theorem, for $H<G$ a finite subgroup, every
$y\in\mathbb{Z}[1/\left\vert H\right\vert ]\otimes R_{\mathbb{C}}(H)$ is a
$\mathbb{Z}[1/\left\vert H\right\vert ]$-linear combination of images
$f_{\ast}y_{\alpha}$ with $y_{\alpha}\in\mathbb{Z}[1/\left\vert H\right\vert
]\otimes R_{\mathbb{C}}(C_{\alpha})$, and $C_{\alpha}<H$ a cyclic subgroup.
This means that for $H\in\mathfrak{F}$ one has:
\begin{align*}
&  \left.  H_{\ast}(X^{H}/C_{G}(H);\Lambda_{G})\otimes\left(  \Lambda
_{G}\otimes R_{\mathbb{C}}(H)\right)  \right/  \sim\\
&  =\sum\left.  H_{\ast}(X^{C}/C_{G}(C);\Lambda_{G})\otimes\left(  \Lambda
_{G}\otimes R_{\mathbb{C}}(C)\right)  \right/  \sim
\end{align*}
where the sum is taken over finite cyclic subgroups of $H$, and thus%
\[
H_{2i}^{\mathfrak{Sub}}(X;\,\Lambda_{G}\otimes R_{\mathbb{C}})=0,\quad
\text{for all }i>0
\]
since according to our assumption, the groups $H_{2i}(X^{C}/C_{G}%
(C);\Lambda_{G})$ are zero for all $i>0$. This implies that the
Atiyah-Hirzebruch spectral sequence for $X$ collapses when the coefficients
are tensored with $\Lambda_{G}$, since all differentials either originate or
end up in $\{0\}$. In particular, the edge homomorphism
\[
c_{0}(X):H_{0}^{\mathfrak{F}}(X;\,\Lambda_{G}\otimes R_{\mathbb{C}%
})\longrightarrow K_{0}^{G}(X)\otimes\Lambda_{G}%
\]
is an isomorphism. Now, from the definitions and the fact that for a subgroup
$H$ of $G$ we have $(X^{0})^{H}=(X^{H})^{0}$, we readily obtain the
surjectivity of%
\[
\iota_{\ast}:H_{0}^{\mathfrak{F}}(X^{0};\,\Lambda_{G}\otimes R_{\mathbb{C}%
})\longrightarrow H_{0}^{\mathfrak{F}}(X;\,\Lambda_{G}\otimes R_{\mathbb{C}})
\]
from its counterpart in ordinary homology. Hence the commutative diagram
expressing naturality%
\[
\begin{CD}
{H_0^{{\mathfrak F}}(X^0;\Lambda_G\otimes R_\mathbb{C})} @>c_0(X^0)>> {K^G_0(X^0)\otimes\Lambda_G}\\
@V{\iota}_*VV @V{\iota}_*VV\\{H_0^{\mathfrak{F}}(X;\Lambda_G\otimes R_\mathbb{C})}@>c_0(X)>>{K_0^G(X)\otimes\Lambda_G}
\end{CD}\label{c_0 from X^0 to X}%
\]
implies the surjectivity of%
\[
\iota_{\ast}:K_{0}^{G}(X^{0})\otimes\Lambda_{G}\longrightarrow K_{0}%
^{G}(X)\otimes\Lambda_{G}\,,
\]
which establishes that $X$ is $K_{0}^{G}\Lambda_{G}$-discrete.
\end{proof}

\begin{cor}
\label{acyclic=>Ldiscrete} For any pervasively acyclic group $G$,
$\underline{E}G$ is $K_{0}^{G}\Lambda_{G}$-discrete.
\end{cor}


\begin{proof}
For every finite cyclic subgroup $C<G$, $(\underline{E}G)^{C}$ is a model for
$\underline{E}C_{G}(C)$. The natural map $EC_{G}(C)\rightarrow\underline
{E}C_{G}(C)$ gives rise to a map
\[
BC_{G}(C)=\left.  EC_{G}(C)\right/  C_{G}(C)\rightarrow\left.  \underline
{E}C_{G}(C)\right/  C_{G}(C)=(\left.  \underline{E}G)^{C}\right/  C_{G}(C)
\]
which induces an isomorphism
\[
H_{\ast}(BC_{G}(C);\Lambda_{G})\cong H_{\ast}(\left.  (\underline{E}%
G)^{C}\right/  C_{G}(C);\Lambda_{G}),
\]
by Lemma \ref{BrownSpectral} applied to the group $C_{G}(C)$. Thus, since the
centralizers $C_{G}(C)$ are acyclic, the result follows by choosing
$X=\underline{E}G$ in the previous theorem.
\end{proof}


\section{Pervasively acyclic groups}

In this section we introduce a functorial embedding of a given group $G$ in an
acyclic group $A(G)$ that has further strong properties required for our
arguments. Embeddings into acyclic groups have historically been important for
algebraic $K$-theory \cite{Wagoner} and algebraic topology \cite{KanThurston}.
One of our further requirements for $A(G)$ is that the centralizers in $A(G)$
of finitely generated abelian subgroups of $A(G)$ should also be acyclic (so
that $A(G)$ is pervasively acyclic, cf.~Definition \ref{pervasively}). In the
extreme case where $G$ is itself abelian, a weaker form of this was already
known to be possible by making $G$ the centre of an acyclic group \cite{BDH},
\cite{Berr: two functors}, \cite{BerrSGTC}. A prominent class of acyclic
groups suited to our purpose comprises the binate groups \cite{BerrSGTC}. We
now recall the definition.

\begin{defi}
A group $G$ is said to be \emph{binate} if for any finitely generated subgroup
$H$ of $G$ there is a homomorphism $\varphi_{H}:H\rightarrow G$ and an element
$u_{H}\in G$ such that for all $h$ in $H$ we have%
\[
h=[u_{H},\,\varphi_{H}(h)]=u_{H}\varphi_{H}(h)u_{H}^{-1}\varphi_{H}(h)^{-1}.
\]

\end{defi}

Obviously $\varphi_{H}$ is injective, while the fact that it is a homomorphism
implies, from the usual product formula for commutators, that its image
commutes with $H.$ So this apparatus embeds a pair of commuting copies of each
finitely generated subgroup of $G$ in $G$ (\emph{binate }= arranged in pairs).
The key property of binate groups is their acyclicity.

\begin{thm}
[see \cite{Berrick: bk}, (11.11)]\label{binate is acyclic}Every binate group
is acyclic.
\end{thm}

In general, there is a construction for embedding a given group in a binate
group, the \emph{universal binate tower} \cite{BerrSGTC}. It is universal in
the sense that it maps to any other such construction \cite{BerrVarad}. See
\cite{Berrick: dichotomy} for justification of the role of binate groups in
this context. The construction we now give is an adaptation of the universal
binate tower.

In the following, we use the notation
\begin{align*}
\Delta_{G}  &  =\left\{  (g,g)\in G\times G\mid g\in G\right\}  ,\\
\Delta_{F}^{^{\prime}}  &  =\left\{  (f^{-1},f)\in F\times F\mid f\in
F\right\}  .
\end{align*}

\begin{defi}
Let $H$ be a group with $\{F_{i}\}_{i\in I}$ as the set of all its finitely
generated abelian subgroups. For $i\in I$, write $C_{i}$ for the centralizer
in $H$ of $F_{i}$. Then in $H\times H$ the subgroups $(1\times F_{i}%
)=\{(1,f)\mid f\in F_{i}\}$ and $\Delta_{C_{i}}$ commute, so that their
product $(1\times F_{i})\cdot\Delta_{C_{i}}$ is also a subgroup of $H\times
H$. Likewise, $\Delta_{F_{i}}^{\prime}\cdot(1\times C_{i})$ is also a
subgroup, and the obvious bijection
\begin{align*}
(1\times F_{i})\cdot\Delta_{C_{i}}  &  \longrightarrow\Delta_{F_{i}}^{\prime
}\cdot(1\times C_{i})\\
(k,fk)  &  \longmapsto(f^{-1},fk)
\end{align*}
is a group isomorphism. Now define $A_{1}(H)$ to be the generalized HNN
extension
\[
A_{1}(H)=\mathrm{HNN}(H\times H;\ (1\times F_{i})\cdot\Delta_{C_{i}}%
\cong\Delta_{F_{i}}^{\prime}\cdot(1\times C_{i}),\ t_{i})_{i\in I}%
\]
meaning that, whenever $f\in F_{i}$ and $k\in C_{i},$
\[
(k,fk)=t_{i}(f^{-1},fk)t_{i}^{-1}.
\]

\end{defi}


\begin{lemma}
\label{01o05a}The inclusion $h\mapsto(h,1)$ of $H$ in $H\times H$ as
$H\times1$ extends to a functorial inclusion of $H$ in $A_{1}(H)$.
\end{lemma}

\begin{proof}
That we obtain an inclusion is an application of the Higman-Neumann-Neumann
Embedding Theorem for HNN\ extensions to the present situation. Functoriality
is a routine check, since homomorphisms map finitely generated abelian
subgroups to finitely generated abelian subgroups and centralizers into centralizers.
\end{proof}

The inclusion of $H$ in $A_{1}(H)$ is used implicitly in the sequel. Note that
an element in $A_{1}(H)$ can be written as a word involving elements of
$H\times H$ and stable letters.

\begin{lemma}
\label{01o05b}For all $i$ in $I$, $t_{i}$ centralizes $F_{i}$, that is,
$t_{i}\in C_{A_{1}(H)}(F_{i})$.
\end{lemma}

\begin{proof}
From the construction we have, for any $f\in F_{i}$,
\[
t_{i}(f,1)t_{i}^{-1}=t_{i}(f,f^{-1}f)t_{i}^{-1}=(f,f^{-1}f)=(f,1).
\]

\end{proof}


\begin{lemma}
\label{01o05c}If two elements of $H$ are conjugate in $A_{1}(H)$, then they
are conjugate in $H$.
\end{lemma}

\begin{proof}
Fix $h\in H$. Our aim is to contradict the assertion that there exist elements
of $H\times H$ that are conjugate to $h$ in $A_{1}(H)$ but not in $H\times H$.
This contradiction gives the required result, for it means that whenever
$(h^{\prime},h^{\prime\prime})$ in $H\times H$ is conjugate to $h=(h,1),$ then
there exists $x=(x^{\prime},x^{\prime\prime})\in H\times H$ such that
\[
(h^{\prime},h^{\prime\prime})=(x^{\prime},x^{\prime\prime})(h,1)(x^{\prime
},x^{\prime\prime})^{-1}.
\]
(Of course, we are primarily interested in the case where we begin with the
further assumption that $h^{\prime\prime}=1;$ however, the more general case
is used later.) The above equation yields that indeed $h^{\prime\prime}=1$ and
so $(h^{\prime},h^{\prime\prime})=h^{\prime}=x^{\prime}h(x^{\prime})^{-1}$,
making $h$ and $xhx^{-1}$ conjugate in $H$ itself. Therefore, in order to
establish a contradiction, let us take $s$ to be the minimal number of indices
$i\in I$ of stable letters $t_{i}$ occurring in any expression of $x$ as a
word in $A_{1}(H)$, where $x$ varies among all elements of $A_{1}(H)$ for
which $xhx^{-1}$ lies in $H\times H$ but fails to be conjugate to $h$ in
$H\times H$. The fact that $x\notin H\times H$ means $s\geq1$. Let $w$ be such
a word involving precisely $s$ distinct stable letters.

Write $\bar{h}=whw^{-1}$, and let $i(1),\ldots,i(s)\in I$ be the indices of
stable letters $t_{i}$ occurring in $w$. We argue as in Step 2 of the proof of
Britton's Lemma \cite{Britton}. We first observe that the trivial word
$\bar{h}^{-1}whw^{-1}$ lies in the subgroup
\[
B_{s}=\mathrm{HNN}(H\times H;\ (1\times F_{i(j)})\cdot\Delta_{C_{i(j)}}%
\cong\Delta_{F_{i(j)}}^{\prime}\cdot(1\times C_{i(j)}),\ t_{i(j)}%
)_{j=1,\ldots,s}%
\]
of $A_{1}(H)$. We put $B_{0}=H\times H,$ and, for $r=1,\ldots,s$, write
\[
B_{r}=\mathrm{HNN}(B_{r-1};\ (1\times F_{i(r)})\cdot\Delta_{C_{i(r)}}%
\cong\Delta_{F_{i(r)}}^{\prime}\cdot(1\times C_{i(r)}),\ t_{i(r)}).
\]
Then, once again using the Embedding Theorem, we have
\[
H=H\times1<H\times H=B_{0}<B_{1}<\cdots<B_{r-1}<B_{r}<\cdots<B_{s}.
\]
Clearly it suffices to establish the following.

\bigskip\noindent\textbf{Claim. }\emph{If }$h$ \emph{is conjugate in }$B_{r}$
$(1\leq r\leq s)$ \emph{to an element }$g$ \emph{of }$B_{r-1},$ \emph{then it
is conjugate in }$B_{r-1}$ \emph{to }$g.$ \bigskip

To prove this, we write $F=F_{i(r)}$, $C=C_{i(r)},\ t=t_{i(r)}$. By the Normal
Form Theorem \cite{LyndonSchupp} p.182 for $B_{r}$ as an HNN extension of
$B_{r-1}$, we may suppose that $n$ is minimal among all $x\in B_{r}$ with
$xhx^{-1}=g$ and
\[
x=g_{0}t^{\varepsilon_{1}}g_{1}t^{\varepsilon_{2}}g_{2}\cdots t^{\varepsilon
_{n}}g_{n}%
\]
reduced. That is, each $g_{i}\in B_{r-1}$, each $\varepsilon_{i}\in\{\pm1\}$,
and there is no subword of the form
\[
tg_{i}t^{-1}\quad\text{with\quad}g_{i}\in\Delta_{F}^{\prime}\cdot(1\times C)
\]
nor
\[
t^{-1}g_{i}t\quad\text{with\quad}g_{i}\in(1\times F)\cdot\Delta_{C}.
\]
To show that in fact $x\in B_{r-1}$ we contradict the minimality of $n$ when
$n\geq1$. Since in $B_{r}$
\[
1=g^{-1}g_{0}t^{\varepsilon_{1}}g_{1}\cdots t^{\varepsilon_{n}}g_{n}%
hg_{n}^{-1}t^{-\varepsilon_{n}}\cdots g_{1}^{-1}t^{-\varepsilon_{1}}g_{0}%
^{-1},
\]
by Britton's Lemma (see \cite{LyndonSchupp}) we must have either
\[
\varepsilon_{n}=1\quad\text{and\quad}g_{n}hg_{n}^{-1}\in\Delta_{F}^{\prime
}\cdot(1\times C)
\]
or
\[
\varepsilon_{n}=-1\quad\text{and\quad}g_{n}hg_{n}^{-1}\in(1\times
F)\cdot\Delta_{C}.
\]
Since $g_{n}\in B_{r-1}\leq B_{s-1}$ and $g_{n}hg_{n}^{-1}\in H\times H$, by
minimality of $s$ we must have
\[
g_{n}hg_{n}^{-1}=yhy^{-1}\quad\text{with\quad}y\in H\times H.
\]
As noted at the beginning, this forces $g_{n}hg_{n}^{-1}\in H$. Therefore
\[
g_{n}hg_{n}^{-1}\in H\mathbb{\cap}((\Delta_{F}^{\prime}\cdot(1\times
C))\mathbb{\cup}((1\times F)\cdot\Delta_{C}))=F.
\]
However, by Lemma \ref{01o05b}, $t\in C_{B_{r}}(F)$. Hence
\begin{align*}
xhx^{-1}  &  =g_{0}t^{\varepsilon_{1}}g_{1}t^{\varepsilon_{2}}g_{2}\cdots
t^{\varepsilon_{n-1}}g_{n-1}g_{n}hg_{n}^{-1}g_{n-1}^{-1}t^{-\varepsilon_{n-1}%
}\cdots g_{1}^{-1}t^{-\varepsilon_{1}}g_{0}^{-1}\\
&  =zhz^{-1}%
\end{align*}
where
$z=g_{0}t^{\varepsilon_{1}}g_{1}t^{\varepsilon_{2}}g_{2}\cdots
t^{\varepsilon_{n-1}}g_{n-1}g_{n}$ involves fewer than $n$
occurrences of $t$, thereby contradicting the minimality of $n$.
\end{proof}

The final property that we need is that no new primary torsion is created in
the construction.

\begin{lemma}
\label{01o05d}$A_{1}(H)$ contains an element of prime power order $p^{k}$, if
and only if $H$ also contains an element of order $p^{k}$.
\end{lemma}


\begin{proof}
It is evident that $H$ and $H\times H$ have the same prime powers arising as
orders of elements. The Torsion Theorem for HNN extensions \cite{LyndonSchupp}
p.185, after iteration as in the proof of Lemma \ref{01o05c}, shows that
$H\times H$ and $A_{1}(H)$ share the same finite orders of elements. (Of
course $H\times H$, and hence $A_{1}(H)$, may have non-prime-power orders not
found in $H$ itself.)
\end{proof}


\begin{defi}
\label{Strongly acyclic group}Let $G$ be a group. We write $A_{1}=A_{1}(G),$
and for $n\geq2$ inductively define $A_{n}(G)=A_{1}(A_{n-1})$, which we also
write as $A_{n}$. Thus by Lemma \ref{01o05a} $A_{n-1}\leq A_{n}$, and we put
$A=A(G)=\mathbb{\cup}A_{n}$.
\end{defi}


\begin{thm}
\label{Berrick} The homomorphism $G\rightarrow A(G)$ has the following properties.

\begin{enumerate}
\item[(a)] It is a functorial inclusion.

\item[(b)] Every finitely generated abelian subgroup of $A(G)$ has its
centralizer in $A(G)$ binate, hence acyclic. In particular, $A(G)$ is
pervasively acyclic.

\item[(c)] If two elements of $G$ are conjugate in $A(G),$ then they are
conjugate in $G$ itself.

\item[(d)] The prime powers that occur as orders of elements of $A(G)$ are
precisely those that occur as orders of elements of $G$.

\item[(e)] If $G$ is countable, $A(G)$ is too.
\end{enumerate}
\end{thm}


\begin{proof}
(a), (c) and (d) are easy consequences of Lemmas \ref{01o05a}, \ref{01o05c}
and \ref{01o05d} respectively. So we concentrate on proving (b). In view of
Theorem \ref{binate is acyclic}, we show that, for each finitely generated
abelian subgroup $F$ of $A(G)$, the centralizer $C_{A}(F)$ is binate.
Therefore let $H$ be a finitely generated subgroup of $C_{A}(F)$. It follows
that we can find $n$ with both $H$ and $F$ subgroups of $A_{n}$, and so that
$F=F_{i}$ for some suitable member $i$ of the index set of finitely generated
abelian subgroups of $A_{n}$. Then the homomorphism $\varphi$ sending $h\in H$
to $(1,h)\in A_{n+1}$ maps into $C_{A}(F)$. By Lemma \ref{01o05b} we have also
$t_{i}\in C_{A}(F)$. Thus the fact that, for any $h\in H$,
\[
h=(h,h)(1,h)^{-1}=t_{i}(1,h)t_{i}^{-1}(1,h)^{-1}=[t_{i},\,\varphi(h)]
\]
reveals $C_{A}(F)$ to be binate, as required.

For (e) one notices that a countable group has only countably many finitely
generated subgroups.
\end{proof}


\section{The Bass conjecture for amenable groups: proof of Theorem
\ref{Bost=>L1Bass copy(1)}}

In this section we show how the use of pervasively acyclic groups allows us to
deduce the Bass conjecture from the Bost conjecture on the assembly map
\[
\beta_{\ast}^{G}:K_{\ast}^{G}(\underline{E}G)\rightarrow K_{\ast
}^{\mathrm{top}}(\ell^{1}(G)).
\]
Notice that in the case where $\ast=0$ and $A$ is any Banach algebra, the
groups $K_{0}(A)$ and $K_{0}^{\mathrm{top}}(A)$ agree, as they both are
defined to be the Grothendieck group of finitely generated projective
$A$-modules -- the topology does not enter in this case. \medskip

\noindent\textit{Proof of Theorem \ref{Bost=>L1Bass copy(1)}. }Consider the
embedding $G\rightarrow A$ where $A=A(G)$ is the pervasively acyclic group of
Definition \ref{Strongly acyclic group}. This embedding, as well as the
inclusion $\iota:\underline{E}A^{0}\rightarrow\underline{E}A$, yields the
following commutative diagram:
\[
\begin{CD}
K_0^G(\underline{E}G) @>{\beta_0^G}>> K_0(\ell^1(G)) @>{HS^1}>> {\ell^1([G])}\\
@VVV                                   @VVV                  @VVV  \\
K_0^A(\underline{E}A) @>{\beta_0^A}>> K_0(\ell^1(A)) @>{HS^1}>> {\ell^1([A])}\\
@A{\iota_*}AA                          @AAA                  @A{j}AA  \\
K_0^A(\underline{E}A^0) @>{\beta}>> \bigoplus R_{\mathbb{C}}F_\alpha @>{\gamma}>>\bigoplus\ell^1([F_\alpha]).\\
\end{CD}\label{beta and HS^1 from G to A}%
\]
The map $\beta$ arises from Remark \ref{0-sk}, by writing $\underline{E}A^{0}$
as a disjoint union of orbits $A/F_{\alpha}$ (where the $F_{\alpha}$ are the
finite subgroups of $A$). We now discuss the maps $\gamma$ and $j$. For
$F_{\alpha}<A$ finite, we have a commutative diagram
\[
\begin{CD}
K_0(\ell^1(A))           @>{HS^1}>>                  {\ell^1([A])}\\
@AAA                                            @A{j_\alpha}AA\\
R_\mathbb{C}(F_\alpha)=K_0(\ell^1(F_\alpha))@>{HS^1}>>\ell^1([F_\alpha])\\
\end{CD}\label{HS^1 from F_alpha to A}%
\]
Here $j_{a}$ is induced by $[F_{\alpha}]\rightarrow\lbrack A]$, and therefore
maps into $\bigoplus_{\mathrm{FC}(A)}\mathbb{C}$, where $\mathrm{FC}(A)$
denotes the set of conjugacy classes of finite order elements of $A$. The map
$\gamma$ is just the sum of Hattori-Stallings traces (on the $F_{\alpha}$),
and $j:\bigoplus\ell^{1}([F_{\alpha}])\rightarrow\ell^{1}([A])$ is the map
that restricts on each $\ell^{1}([F_{\alpha}])$ to $j_{\alpha}$.

\smallskip

After tensoring the left two columns of this diagram with $\mathbb{Q}$, by
assumption the Bost assembly map $\beta_{0}^{G}\otimes\mathrm{Id}_{\mathbb{Q}%
}$ becomes surjective, and the map $\iota_{\ast}\otimes\mathrm{Id}%
_{\mathbb{Q}}$ becomes an epimorphism because of Corollary
\ref{acyclic=>Ldiscrete}. It is a consequence of Theorem \ref{Berrick} part
(c) that $[G]\subseteq\lbrack A]$, and thus
\[
\ell^{1}([G])\subset\ell^{1}([A])\supset\ell^{1}\left(  FC(A)\right)
\supset\bigoplus\limits_{\mathrm{FC}(A)}\mathbb{C},
\]
with
\[
\bigoplus\limits_{\mathrm{FC}(G)}\mathbb{C}=\ell^{1}([G])\bigcap
\bigoplus\limits_{\mathrm{FC}(A)}\mathbb{C}.
\]
Therefore the map
\[
{HS^{1}}:K_{0}(\ell^{1}(G))\rightarrow{\ell^{1}([G])}%
\]
takes its values in the $\mathbb{C}$-vector space spanned by the conjugacy
classes of elements of finite order.\hfill$\square$\medskip

\section{On the image of the Kaplansky trace}

\label{ImageTrace}
We now turn to the idempotent conjecture for $\ell^{1}(G)$ and the proof of
Corollary \ref{B}. Using similar arguments to those that have been used in part 1 of
\cite{Mislin}, Lemma 7.2, one could deduce the idempotent conjecture for
$\ell^{1}(G)$ from the surjectivity of the Bost assembly map. However, our aim
is to prove something stronger, namely the idempotent conjecture for $\ell
^{1}(G)$ out of the \emph{rational surjectivity} of the Bost assembly map. We
start with some standard facts concerning the Kaplansky trace.

\smallbreak

Recall that $C_{r}^{\ast}(G)$ acts on $\ell^{2}(G)$. Considering $1\in G$ as
an element of $\ell^{2}(G)$, one defines for $f\in C_{r}^{\ast}(G)$ the
\emph{Kaplansky trace}
\[
\kappa(f)=\left\langle f(1),1\right\rangle \in\mathbb{C}.
\]
This defines a trace $\kappa:C_{r}^{\ast}(G)\rightarrow\mathbb{C}$ which
extends the Kaplansky trace $\ell^{1}(G)\rightarrow\mathbb{C}$ described
earlier, in the sense that the resulting map in $K$-theory (which we still
denote by $\kappa$) fits into a commutative diagram
\[
\begin{CD}
{K_0(\ell^1(G))}    @>\kappa>>       \mathbb{C}     \\
@VVV                            @|\\
{K_0(C^*_r(G))}        @>\kappa>>           \mathbb{C}.    \\
\end{CD}\label{kappa from l^1 to C^*_r}%
\]
Moreover, $\kappa$ is natural in the sense that for $H<G$ one has a
commutative diagram
\[
\begin{CD}
{K_0(C^*_r(H))}    @>\kappa>>       \mathbb{C}     \\
@VVV                            @|\\
{K_0(C^*_r(G))}        @>\kappa>>           \mathbb{C}.    \\
\end{CD}\label{kappa from H to G}%
\]
Combined with Theorem \ref{Bost=>L1Bass copy(1)}, the following proposition
completes the proof of Corollary \ref{B}.

\begin{prop}
Let $G$ be a torsion-free group satisfying Conjecture~\ref{L1Bass}. Then
$\ell^{1}(G)$ contains no idempotent other than $0$ and $1$.
\end{prop}

\begin{proof}
We first show that in the case where Conjecture \ref{L1Bass} holds and $G$ is
torsion-free, both the Kaplansky trace $\kappa$ and the augmentation trace
$\epsilon$, considered as maps $K_{0}(\ell^{1}(G))\to\mathbb{C}$, coincide. On
the one hand, for $G$ torsion-free, Conjecture \ref{L1Bass} says that
$HS^{1}=\epsilon_{[1]}=\epsilon$ (since $G$ being torsion-free implies that
$\mathrm{FC}(G)=\{1\}$), and on the other hand as observed earlier one has
that $\kappa=\epsilon_{[1]}$. Since the augmentation trace $\epsilon
:K_{0}(\ell^{1}(G))\to\mathbb{C}$ is known to assume integral values only
(recall from Section 2 that it factors through $K_{0}(\mathbb{C})$), this
means that the Kaplansky trace $\kappa$ has integral range as well.

If $a\in\ell^{1}(G)\subset C_{r}^{\ast}(G)$ is an idempotent, then the
projective ideal $P=C_{r}^{\ast}(G){\cdot}a\subset C_{r}^{\ast}(G)$ satisfies
$\kappa([P])=0$ or $1$, so that $P=0$ or $C_{r}^{\ast}(G)$. Therefore $a=0$ or
$1$ by the usual argument concerning idempotents in $C_{r}^{\ast}(G)$ (cf.
\cite{BCH}, Proposition 7.16).
\end{proof}

Using similar ideas to those in the previous section, we now give another
application of pervasively acyclic groups and Proposition
\ref{acyclic=>Ldiscrete}, but in the context of the Baum-Connes conjecture.
Here we recall:

\begin{con}
[Baum-Connes]\label{BC}Let $G$ be a countable discrete group. The Baum-Connes
assembly map
\[
\mu^{G}_{*}:K_{*}^{G}(\underline{E}G)\to K^{\mathrm{top}}_{*}(C^{*}_{r}(G))
\]
is an isomorphism.
\end{con}

We refer to Baum, Connes and Higson \cite{BCH}, as well as
Lafforgue's work \cite{Lafforgue} for the construction of the
Baum-Connes assembly map, and for considerable partial information
concerning the validity of the Baum-Connes conjecture. Our
techniques allow us to give a new proof of a recent theorem of
L\"uck \cite{Lueck}.

\begin{thm}
\label{Lk}Let $\kappa: K_{0}(C^{*}_{r}(G))\rightarrow \mathbb{C}$
be the Kaplansky trace and $\mu^{G}_{0}:
K_{0}^{G}(\underline{E}G)\rightarrow K_{0}(C^{*}_{r}(G))$ the Baum
Connes assembly map. Then the image of $\kappa\circ\mu^{G}_{0}$ is
contained in $\Lambda_{G}$.
\end{thm}

\begin{proof}
Consider the embedding $G\rightarrow A$ where $A=A(G)$ denotes the pervasively
acyclic group of Definition \ref{Strongly acyclic group}. Notice that
$\Lambda_{G}=\Lambda_{A(G)}$ by Theorem \ref{Berrick} (d). This equality, as
well as the inclusion $\iota:\underline{E}A^{0}\rightarrow\underline{E}A$,
yields the commutative diagram
\[
\begin{CD}
{K_0^G(\underline{E}G)} @>{\mu^G}>> {K_0(C^*_r(G))}      @>\kappa>>          \mathbb{C}     \\
@VVV                                 @VVV                            @V{\rm Id}VV\\
{K_0^A(\underline{E}A)} @>{\mu^A}>> {K_0(C^*_r(A))}      @>\kappa>>           \mathbb{C}     \\
@A{\iota_*}AA                        @AAA                             @AAA       \\
{K_0^A(\underline{E}A^0)} @>{\beta}>> {\bigoplus R_{\mathbb{C}}F_\alpha} @>{\sigma}>> \Lambda_G\\
\end{CD}\label{mu and kappa from G to A}%
\]
where $\beta$ is as in the proof of Theorem \ref{Bost=>L1Bass copy(1)}, and
$\sigma$ is given as follows. If $\{x_{\alpha}\}$ is an element in ${\bigoplus
R_{\mathbb{C}}F_{\alpha}}$, then
\[
\sigma(\{x_{\alpha}\})=\sum_{\alpha}\kappa_{\alpha}(x_{\alpha})
\]
where $\kappa_{\alpha}:R_{\mathbb{C}}F_{\alpha}\rightarrow\mathbb{C}$ is the
Kaplansky trace on $K_{0}(C_{r}^{\ast}(F_{\alpha}))=R_{\mathbb{C}}F_{\alpha}$,
which takes its values in $\frac{1}{|F_{\alpha}|}\mathbb{Z}\subset\Lambda_{G}%
$. Indeed, the groups $F_{\alpha}$ are finite, and for a finitely generated
projective $\mathbb{C}F_{\alpha}$-module $P$ one has, according to Bass
(\cite{bass}, Corollary 6.3), by comparing the Kaplansky trace of $P$ with the
Kaplansky trace of $P$ considered as a module over $\mathbb{C}\{1\}$,
\[
\kappa([P])=\frac{1}{|F_{\alpha}|}\mathrm{dim}_{\mathbb{C}}P\in\frac
{1}{|F_{\alpha}|}\mathbb{Z}.
\]
Now write $M$ for the image of $\kappa\circ\mu^{G}_{0}$. As the
map $\iota_{\ast} \otimes\mathrm{Id}_{\Lambda_{G}}$ becomes an
epimorphism because of Corollary \ref{acyclic=>Ldiscrete},
commutativity of the diagram now shows that $M\otimes\Lambda_{G}$
lies in $\Lambda_{G}\otimes\Lambda_{G}$. (This is essentially
L\"uck's theorem.) The $\mathbb{Z}$-flatness of submodules of
$\mathbb{C}$ then gives injections
\[
M\cong{M\otimes\mathbb{Z}}\rightarrow{M\otimes\Lambda_{G}}
\rightarrow\Lambda_{G}\otimes\Lambda_{G}\cong\Lambda_{G},
\]
and so the result.
\end{proof}


\section{New examples of groups satisfying the Bass conjecture}

\label{ClasseC'}We describe a wide class of groups for which the Bost
conjecture is known. This class, which Lafforgue in \cite{Lafforgue} called
$\mathcal{C}^{\prime}$, includes all discrete countable groups acting
metrically properly and isometrically on one of the following spaces:

\begin{itemize}
\item[(a)] an affine Hilbert space (those groups are said to have the
\emph{Haager\-up property}, or to be \emph{a-T-menable});

\item[(b)] a uniformly locally finite, weakly $\delta$-geodesic and strongly
$\delta$-bolic space (we will see that cocompact CAT(0)-groups satisfy this assumption);

\item[(c)] a non-positively curved Riemannian manifold, with curvature bounded
from below and bounded derivative of the curvature tensor (with respect to the
connection induced from the Levi-Civita connection on the tangent bundle).
\end{itemize}

\begin{thm}
\emph{(Lafforgue \cite{Lafforgue})} The Bost conjecture holds for
any group in the class $\mathcal{C}^{\prime}.$
\end{thm}

We now discuss in turn the three classes of groups specified by (a), (b) and
(c) above.

\medskip

\noindent\textbf{Class (a).} Here, we have Theorem \ref{Lafforgue: Bost} as a
special case. We recall that the class of groups satisfying the Haagerup
property has the following closure properties (see \cite{welches}).

\begin{itemize}
\item[--] The Haagerup property is closed under taking subgroups, and direct
products.

\item[--] If $G$ acts on a locally finite tree with finite edge stabilizers,
and with the vertex stabilizers having the Haagerup property, then so does
$G$.

\item[--] If $G=\bigcup_{n\geq0}G_{n}$, with $G_{i}<G_{i+1}$ for all $i$, and
each $G_{i}$ has the Haagerup property, then so does $G$.

\item[--] If $G$ and $H$ are countable amenable groups and $C$ is central in
$H$ and $G$ then $G*_{C}H$ has the Haagerup property (in particular, free
products of countable amenable groups have the Haagerup property).
\end{itemize}

\smallskip

\noindent\textbf{Class (b). }We now turn to the second class of groups
contained in $\mathcal{C}^{\prime}.$

\begin{defi}
A metric space $(X,d)$ is called \emph{uniformly locally finite} if, for any
$r\geq0$, there exists $k\in\mathbb{N}$ such that any ball of radius $r$
contains at most $k$ points (notice that this forces $X$ to be discrete). For
$\delta>0$, the metric space $X$ is termed \emph{weakly $\delta$-geodesic} if,
for any $x,y\in X$ and $t\in[0,d(x,y)]$, there exists $a\in X$ such that%
\[
d(a,x)\leq t+\delta,\ \ d(a,y)\leq d(x,y)-t+\delta.
\]

\end{defi}

The two conditions above are automatically satisfied by any finitely generated
group $G$ endowed with the word metric associated to any finite generating
set, and by the orbit of a point in a Riemannian manifold with non-positive
curvature, under a group acting properly, isometrically and cocompactly on the
manifold. The following definitions are taken from \cite{KS} and
\cite{Lafforgue}.

\begin{defi}
Given $\delta>0$, a metric space $(X,d)$ is said to be \emph{weakly $\delta
$-bolic} if the following conditions are satisfied:

\begin{itemize}
\item[$(b_{1})$] For any $r>0$, there exists $R=R(\delta,r)>0$ such that, for
any four points $x_{1},x_{2},y_{1},y_{2}$ in $X$ satisfying $d(x_{1}%
,y_{1})+d(x_{2},y_{2})\leq r$ and $d(x_{1},y_{2})+d(y_{1},x_{2})\geq R$, one
has:
\[
d(x_{1},x_{2})+d(y_{1},y_{2})\leq d(x_{1},y_{2})+d(y_{1},x_{2})+\delta.
\]

\item[$(b_{2})$] There exists a map $m:X\times X\to X$ such that:

\begin{itemize}
\item[(i)] For all $x,y\in X$,
\[
d(m(x,y),x)\leq\frac{d(x,y)}{2}+\delta\quad\text{and\quad}d(m(x,y),y)\leq
\frac{d(x,y)}{2}+\delta.
\]

\item[(ii)] For all $x,y,z\in X$,
\[
d(m(x,y),z)\leq\max\{d(x,z),d(y,z)\}+2\delta.
\]

\item[(iii)] For all $p\geq0$, there exists $N=N(p)\geq0$ such that, for any
$n\geq N$ and $x,y,z\in X$ with $d(x,z)\leq n$, $d(y,z)\leq n$ and $d(x,y)>n$,
one has
\[
d(m(x,y),z)<n-p.
\]

\end{itemize}
\end{itemize}

Given $\delta>0$, a metric space $(X,d)$ is called \emph{strongly $\delta
$-bolic} if it is weakly $\delta$-bolic and if condition $(b_{1})$ is
satisfied for any $\delta>0$.
\end{defi}

Mineyev and Yu in \cite{MiYu} showed that a hyperbolic group $G$ can be
endowed with a left $G$-invariant metric such that there exists $\delta>0$ for
which $G$ is strongly $\delta$-bolic.

We recall that in a geodesic metric space $(X,d)$ a \emph{geodesic triangle}
$\Delta$ consists of three points $a,b,c$ with three (possibly non-unique)
geodesics joining them, and a \emph{comparison triangle} $\overline{\Delta}$
for $\Delta$ is a euclidean triangle with side lengths $d(a,b),d(b,c),d(c,a)$.
We write $\overline{a},\overline{b}$ and $\overline{c}$ for the vertices of
$\overline{\Delta}$, and if $x$ is a point in $\Delta$ (say on a geodesic
between $a$ and $b$), we write $\overline{x}$ for a comparison point for $x$,
namely a point in $\overline{\Delta}$ such that $d_{E}(\overline{x}%
,\overline{a})=d(x,a)$ (where $d_{E}$ denotes the euclidean distance in
$\mathbb{R}^{2}$). A geodesic metric space $(X,d)$ is termed \textrm{CAT(0)}
if for all geodesic triangles $\Delta$ in X and all $x,y\in\Delta$,
\[
d(x,y)\leq d_{E}(\overline{x},\overline{y})
\]
where $\overline{x},\overline{y}$ are any two comparison points in any
euclidean comparison triangle $\overline{\Delta}$ for $\Delta$. The following is an easy fact.

\begin{prop}
\textrm{CAT(0)} metric spaces are strongly $\delta$-bolic for any $\delta>0$.
\end{prop}

\begin{proof}
We leave to the reader the verification that $\mathbb{R}^{2}$ with its
euclidean distance is a strongly $\delta$-bolic space for any $\delta>0$ (the
rest of the proof relies on this fact, see also \cite{KS} Proposition 2.4.).
Let $(X,d)$ be a \textrm{CAT(0)} metric space. We start by checking that
condition $(b_{2})$ holds. The map $m$ is defined as follows:
\[
m(x,y)=\gamma(t)
\]
where $t=d(x,y)/2$ and $\gamma:[0,d(x,y)]\rightarrow X$ is the
unique geodesic from $x$ to $y$. Point (i) is satisfied by
assumption on $m$, (ii) holds since in \textrm{CAT(0)} spaces the
metric is strictly convex, and (iii) follows from the
\textrm{CAT(0)} inequality. It remains to prove condition
$(b_{1})$ for any $\delta>0$. To do this, we follow Bridson's
advice and use the \textrm{CAT(0)} $4$-points condition (see
Bridson-Haefliger's book \cite{BridsonHafliger} p.164) which says
that in a \textrm{CAT(0)} space, every $4$-tuple of points
$x_{1},x_{2},y_{1},y_{2}$ has a \emph{sub\-embedding} in
$\mathbb{R}^{2}$, meaning that there exist $4$ points
$\overline{x}_{1},\overline{x}_{2},\overline{y}_{1},\overline{y}_{2}$
in
$\mathbb{R}^{2}$ such that $d_{E}(\overline{x}_{i},\overline{y}_{j}%
)=d(x_{i},y_{j})$ for $i=1,2$ and $d(x_{1},x_{2})\leq d_{E}(\overline{x}%
_{1},\overline{x}_{2})$, $d(y_{1},y_{2})\leq d_{E}(\overline{y}_{1}%
,\overline{y}_{2})$. We now take $\delta,r\geq0$ and $R=R(\delta,r)$ as for
$\mathbb{R}^{2}$. By definition of subembedding, the $4$-tuple $\overline
{x}_{1},\overline{x}_{2},\overline{y}_{1},\overline{y}_{2}$ satisfies the
assumptions of $(b_{1})$ as soon as the $4$-tuple $x_{1},x_{2},y_{1},y_{2}$
does, and we conclude with
\begin{align*}
&  d(x_{1},x_{2})+d(y_{1},y_{2})\\
&  \leq d_{E}(\overline{x}_{1},\overline{x}_{2})+d_{E}(\overline{y}%
_{1},\overline{y}_{2})\leq d_{E}(\overline{x}_{1},\overline{y}_{2}%
)+d_{E}(\overline{x}_{2},\overline{y}_{1})+\delta\\
&  =d(x_{1},y_{2})+d(y_{1},x_{2})+\delta.
\end{align*}

\end{proof}

If a countable discrete group $G$ acts properly and cocompactly on a
\textrm{CAT(0)} metric space $X$, then there exists a $\delta>0$ such that,
for any $x_{0}\in X$, $Y=Gx_{0}\subset X$ endowed with the induced metric from
$X$ is weakly $\delta$-geodesic and strongly $\delta$-bolic. Indeed, one may
choose $\delta=2R$, where $R$ is a positive real number (existing by
cocompactness) such that $X=\bigcup\nolimits_{g\in G}B(gx_{0},R).$ Finally,
$G$ is finitely generated (see \cite{BridsonHafliger} p. 439). Thus $Y$ is
automatically uniformly locally finite and, since $G$ obviously acts properly
on $Y$, this shows that cocompact \textrm{CAT(0)} groups are in the class
$\mathcal{C}^{\prime}$.

We recall that that the class of cocompact CAT(0) groups (that is, discrete
groups acting properly, isometrically and cocompactly on a CAT(0) metric
space) is closed under the following operations (see \cite{BridsonHafliger} p.
439):

\begin{itemize}
\item[--] direct products;

\item[--] HNN extensions along finite subgroups;

\item[--] free products with amalgamation along virtually cyclic subgroups.
\end{itemize}

\smallskip

\noindent\textbf{Class (c). }The group
$\mathrm{SL}_{n}(\mathbb{Z})$ (and more generally any discrete
subgroup of a virtually connected semisimple linear Lie group) is
in the class $\mathcal{C}^{\prime}$, as any non-positively curved
symmetric space is a Riemannian manifold satisfying the required
assumptions on the curvature.

\smallskip

\textit{Acknowledgement.} The first-named author thanks the Department of
Mathematics, ETH Z\"{u}rich for its generous and productive hospitality.

\end{document}